\RequirePackage{ifpdf}
\ifpdf
	\documentclass[pdftex]{statsoc}
\else
	\documentclass[dvips]{statsoc}
\fi

\usepackage{amsmath, ifthen}

\newboolean{color-fig}
\setboolean{color-fig}{true} 

\ifpdf
  \usepackage[pdftex]{graphicx}
  \usepackage{subfigure}
  \graphicspath{{pdf-figures/}}
  \usepackage[%
    pdftex,%
    pdftitle = {{Position play in carom billiards as a Markov process}},%
	pdfauthor = {{Mathieu Bouville}},%
	pdfsubject = {{Modeling the 'small games' of carom billiards as Markovian processes in order to account for position play.}},%
	pdfkeywords = {{billiards, sports, Markov process, Bernoulli process}},%
	pdfpagemode = UseOutlines,%
 	breaklinks,%
    hyperindex,%
    bookmarksopen,%
    bookmarksopenlevel=1,%
  ]{hyperref}
\else
  \usepackage[dvips]{graphicx}
  \usepackage{subfigure}
  \usepackage[%
    dvips,%
    pdftitle = {{Position play in carom billiards as a Markov process}},%
	pdfauthor = {{Mathieu Bouville}},%
	pdfsubject = {{Modeling the 'small games' of carom billiards as Markovian processes in order to account for position play.}},%
	pdfkeywords = {{billiards, sports, Markov process, Bernoulli process}},%
	pdfpagemode = UseOutlines,%
 	breaklinks,%
    hyperindex,%
    bookmarksopen,%
    bookmarksopenlevel=1,%
  ]{hyperref}
\fi

\newcommand{\vc}[1]{\mathbf{#1}}
\newcommand{\mt}[1]{\mathit{#1}}

\title[Markovian carom billiards]{\addcontentsline{toc}{chapter}{English}
Position play in carom billiards as a Markov process}%
\author[M. Bouville]{Mathieu Bouville}
\address{Institute of Materials Research and Engineering, 3 Research Link, Singapore 117602\\email address: m-bouville@imre.a-star.edu.sg}
\email{m-bouville@imre.a-star.edu.sg}

\begin{document}

\begin{abstract}
\addcontentsline{toc}{section}{Abstract}
Position play is a key feature of carom billiards: on easy shots players can manage to score while ensuring that the next position will be favorable.  
The difficulty of a shot therefore depends on the previous shot, e.g.\ an easy shot generally follows an easy shot. We introduce a Markov process which accounts for such correlations. This model can explain the long series of easy shots and the high scores which ensue. It also enables us to identify differences in the scoring patterns of players at different skill levels. 
Players can use this model {\itshape via} \url{http://billiards.mathieu.bouville.name/biMar/}.

\end{abstract}

\section{Introduction}
The word ``billiards'' can refer to cue games in general ---including pool and snooker--- or it can refer to a specific game, also called ``carom billiards'' or ``carom''. This game is played with three balls on a pocketless table (a short glossary is provided as an appendix). A player must make his cue ball contact with the two object balls. If he succeeds he scores a point and plays again. At the end of the game an average can be obtained by dividing the number of points by the number of innings. Rankings can be thus established, even for players who never played against one another. This is similar to track and field where times or distances can be compared to rank athletes who never competed against one another.

\subsection{Historical background}
In 1880, in Paris, the American Slosson scored 1~103 points in a row against the French Vignaux who then ran 1~531 points. Ten years later Sch{\"a}fer scored 3~000 points in a row. This corresponds to playing for days without missing. 
Such scores are possible because when players manage to gather the balls close to a rail and keep them together they can have very long runs (series of shots without missing). This technique is called the rail nurse. The players could not play for hours without missing if they had to face many difficult shots. 
To prevent the use of the rail nurse technique ---and thus avoid tediously long matches--- new rules were introduced. In balkline the balls cannot be left in the same area of the table for more than one or two shots; long series (hundreds of points) can still be obtained by playing position. 
In cushion games the cue ball must hit one rail (one-cushion billiards) or three rails (three-cushion billiards) before it contacts the second object ball.

\subsection{\label{uncorrelated_model}Carom billiards as a Bernoulli process}
\begin{figure}
\centering
\setlength{\unitlength}{1cm}
\begin{picture}(14,4)
\subfigure
{
    \label{fig-mu-log(a)}
    \includegraphics[height=4cm]{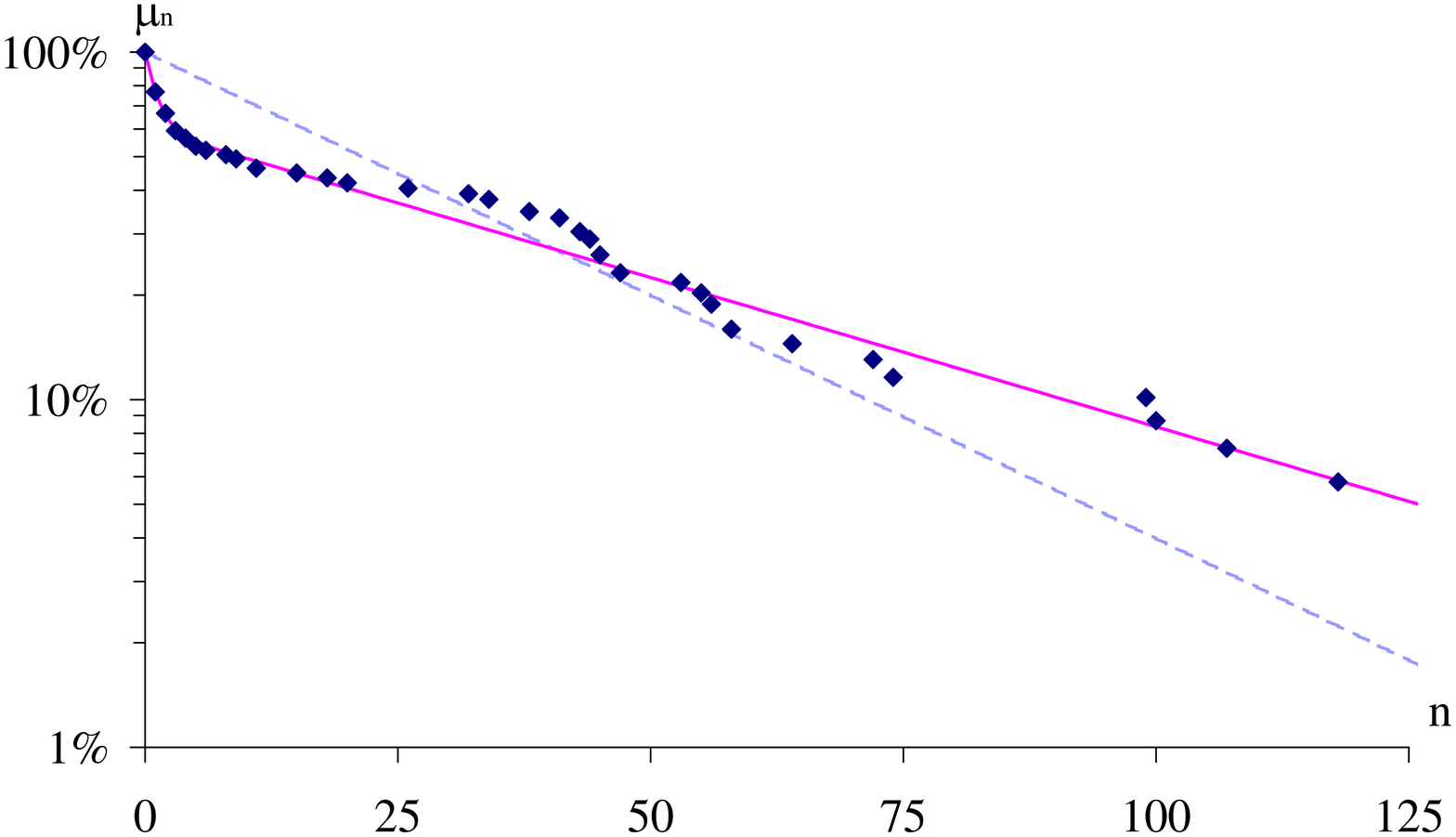}
\put(-0.9,3.5){(a)}
}
\subfigure
{
    \label{fig-mu-log(b)}
    \includegraphics[height=4cm]{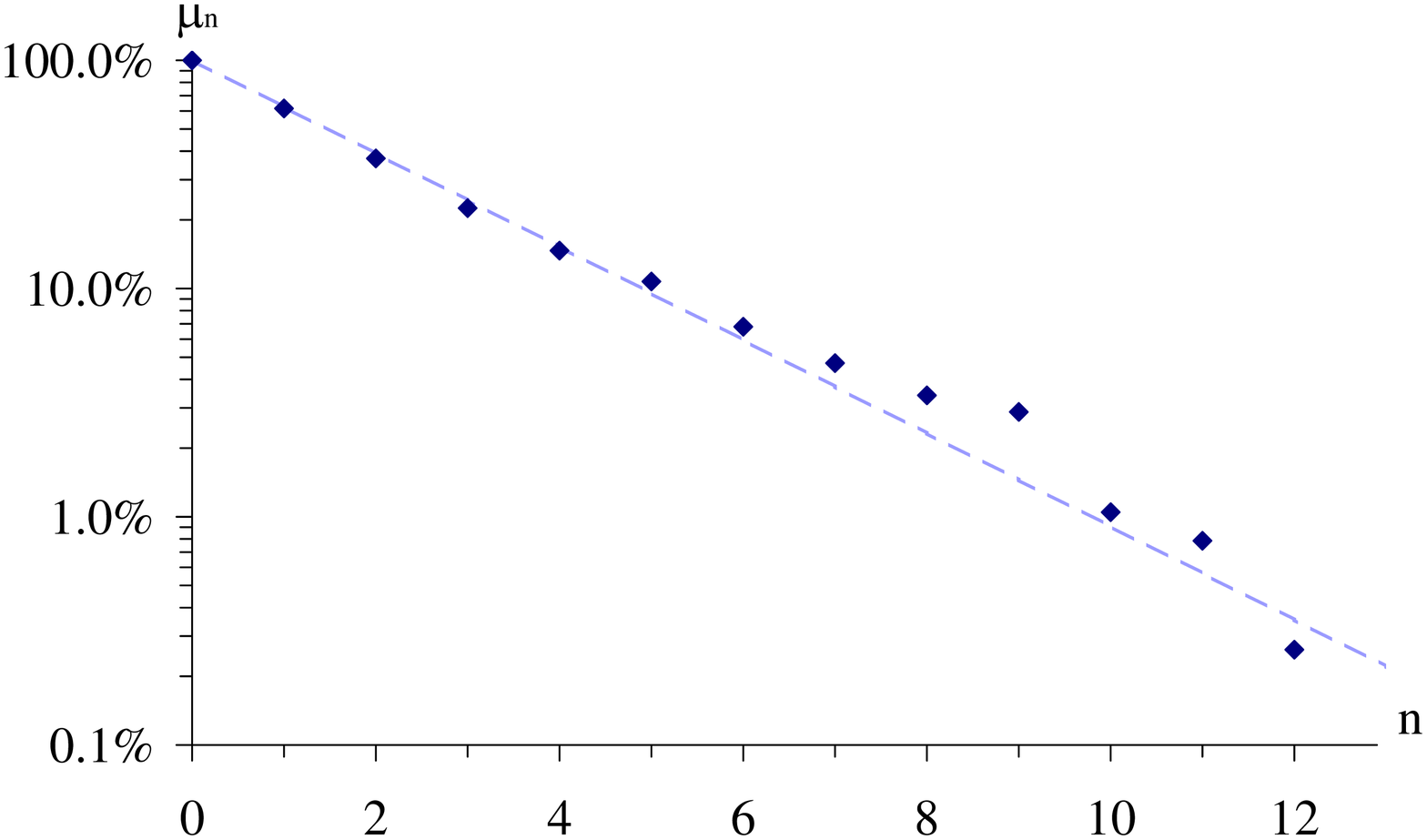}
\put(-0.9,3.5){(b)}
}
\end{picture}
\caption{\label{fig-mu-log}The probability to score at least $n$ points, $\mu_n$, as a function of $n$ for (a) B.~Villiers playing 47/2 balkline and (b) R.~Ceulemans playing three-cushion. Diamonds: data; dashed lines: Bernoullian model; solid line: Markovian model.}
\end{figure}

The game of billiards has been studied in detail from a mechanics point of view by \citet{Coriolis} and more recently by \citet{Petit}. The study of the statistics of billiards on the other hand has been limited to assuming a Bernoulli process. 
 Let $\mu_n$ the probability to score at least $n$ points. In a Bernoulli process with a success rate $\lambda$, for all $n \ge 0$
\begin{equation}
	\mu_n=\lambda^n\!.
	\label{mu_n_Bernoulli}
\end{equation}
\noindent The average length of a run ---subsequently called ``the average''--- is then $m=\lambda/(1-\lambda)$. Figure~\ref{fig-mu-log(a)} shows 47/2 balkline results for Bernard Villiers, based on 69 complete innings (innings which were interrupted because Villiers reached 200 points and won the match are not included). 
The Bernoullian model (dashed line) is a poor fit to the data (diamonds): it overestimates the importance of short runs and underestimates that of long runs.

\subsection{The importance of position play}
A Bernoulli process does not account for position play, which is a very important feature of straight-rail and balkline billiards for instance. The probability to make a shot depends on its difficulty, which is not necessarily random. 
On easy shots, scoring is not an issue and the player can try to play in a way such that the point will be scored and the balls will end in a favorable position for the next shot. This is called playing position. After an easy shot, one can therefore obtain a shot easier than average. 
A player can thus have a series of easy shots and long runs. Such correlations are completely absent from the Bernoullian model described above; it can therefore not properly describe carom billiards.
Even though some limitations of this model have been considered \citep{Fray, Bouville-serie} the effect of position play has never been studied. In order to account for correlations between successive shots we model carom billiards as a Markov process.

\section{Carom billiards as a Markov process}
\subsection{Defining a Markov process}
If one considers that there exist $N_0$ possible positions of the balls on the table, the outcome of a shot can be either success with the balls ending in one of the $N_0$ possible configurations or failure. The probability to score and the positions of the balls for the next shot depend only on the position of the balls at the beginning of the shot, not on where they were earlier in the game.  Thus the succession of shots is a Markov process with $N_0+1$ states
. Unfortunately the transition probabilities are many, $O(N_0{}^2)$, and cannot be readily obtained. 
In order to simplify the model and make it tractable, similar positions can be treated together in order to have only $N$ {\itshape types} of positions. If $N$ is small, the number of parameters is more manageable and these parameters can be calculated and interpreted concretely.

\subsection{The probability to score at least $n$ points}
Let $\vc{P}_{\!n}$ the $N$-dimensional vector whose component $i$ is the probability of having the balls in a position of type~$i$ after shot $n$ (we use ``shot of type~$i$'' and ``position [of the balls on the table] of type~$i$'' as synonyms). 
For all $n\ge0$ the sequence $\left(\vc{P}_{\!n} \right)$ obeys $\vc{P}_{\!n+1} = \mt{K} \vc{P}_{\!n}$.
 The component $k_{ij}$ of the $N \times N$ matrix $\mt{K}$ is the probability to score from a shot of type~$j$ with the next shot being of type~$i$. $\mt{K}$ is \emph{akin} to a transition matrix: the state corresponding to ``misses'' is not included as the corresponding components can straightforwardly be obtained from $\mt{K}$\!. The probability to miss a shot of type~$i$ is $1-\kappa_i$, where $\kappa_i=\sum_j k_{ji}$ is the probability to score on a shot of type~$i$.
 
In what follows we use $N=2$; however the model can straightforwardly be extended to more types of shots. Let $\rho_1 \le \rho_2$ the eigenvalues of $\mt{K}$\!, and $\vc{V}_{\!1}$ and $\vc{V}_{\!2}$ its eigenvectors.
For all $n \ge 0$, $\vc{P}_{\!n}$ can be written as
\begin{equation}
	\vc{P}_{\!n} = \left(\begin{array}{cc}	\vc{V}_{\!1} & \vc{V}_{\!2} \end{array} \right) \left(\begin{array}{cc}
		\rho_1{}^n & 0\\
		 0 & \rho_2{}^n
	\end{array} \right) \left(\begin{array}{cc}	\vc{V}_{\!1} & \vc{V}_{\!2} \end{array} \right)^{-1}
		\left(\begin{array}{cc}	p_0 \\ 1-p_0 \end{array} \right)\negthickspace,
	\label{eq_P}
\end{equation}
\noindent where $p_0$ is the probability that the first shot which the player has to play be difficult. The probability to score at least $n$ points, $\mu_n$, is the sum of the components of $\vc{P}_{\!n}$. For all $n \ge 0$, $\mu_n$ is of the form
\begin{equation}
	\mu_n = (1-Y) \rho_1{}^n + Y \rho_2{}^n\!, 
	\label{mu_Y}
\end{equation}
\noindent where $Y$ is a scalar known as a function of $\vc{V}_{\!1}$, $\vc{V}_{\!2}$, and $p_0$. Figure~\ref{fig-mu-log(a)} shows $\mu_n$ on a semi-logarithmic scale as a function of $n$ for B.~Villiers playing 47/2 balkline: the data follow an exponential law at low $n$ and another exponential at high values of $n$. The Markov process, Eq.~(\ref{mu_Y}), provides a better description of balkline than the Bernoulli process, Eq.~(\ref{mu_n_Bernoulli}). 

The average is the sum over $n$ of $\mu_n$,
\begin{equation}
	m = (1-Y) \frac{\rho_1}{1-\rho_1} + Y \frac{\rho_2}{1-\rho_2}.
	\label{def-m}
\end{equation}
\noindent As $Y$ can be written as a function of $\rho_1$, $\rho_2$, and $m$, Eq.~(\ref{mu_Y}) can be rewritten to have $\mu_n$ as a function of $\rho_1$, $\rho_2$, $m$, and $n$. As the average, $m$, is known from score sheets, $(\mu_n)$ depends on two unknown parameters only, $\rho_1$ and $\rho_2$. We find them by minimizing the sum of the squares of the differences between $\mu_n$ from the data and Eq.~(\ref{mu_Y}).

\subsection{Types of shots}
We considered that there exist two types of shots, defined the matrix $\mt{K}$\!, and used its eigenvectors $\rho_1$ and $\rho_2$ to calculate $\mu_n$ in Eq.~(\ref{mu_Y}). However we never really specified what these two types of shots are. They cannot be random, otherwise the Markov process would boil down to a Bernoulli process. Positions of type~1 (resp.\ 2) must be such that they are generally followed by another position of type~1 (resp.\ 2).
In straight-rail, shots of type~2 would be rail-nurse positions: there is a high probability that a rail-nurse position be followed by another rail-nurse position. Shots of type~1 would be positions of the balls which are less favorable. 
More generally positions of type~1 are ``difficult'' and those of type~2 ``easy''. 

\begin{figure}
\centering
\setlength{\unitlength}{1cm}
\begin{picture}(14,4.6)
\subfigure{
    \label{fig-ratio(a)}
    \includegraphics[height=4.1cm]{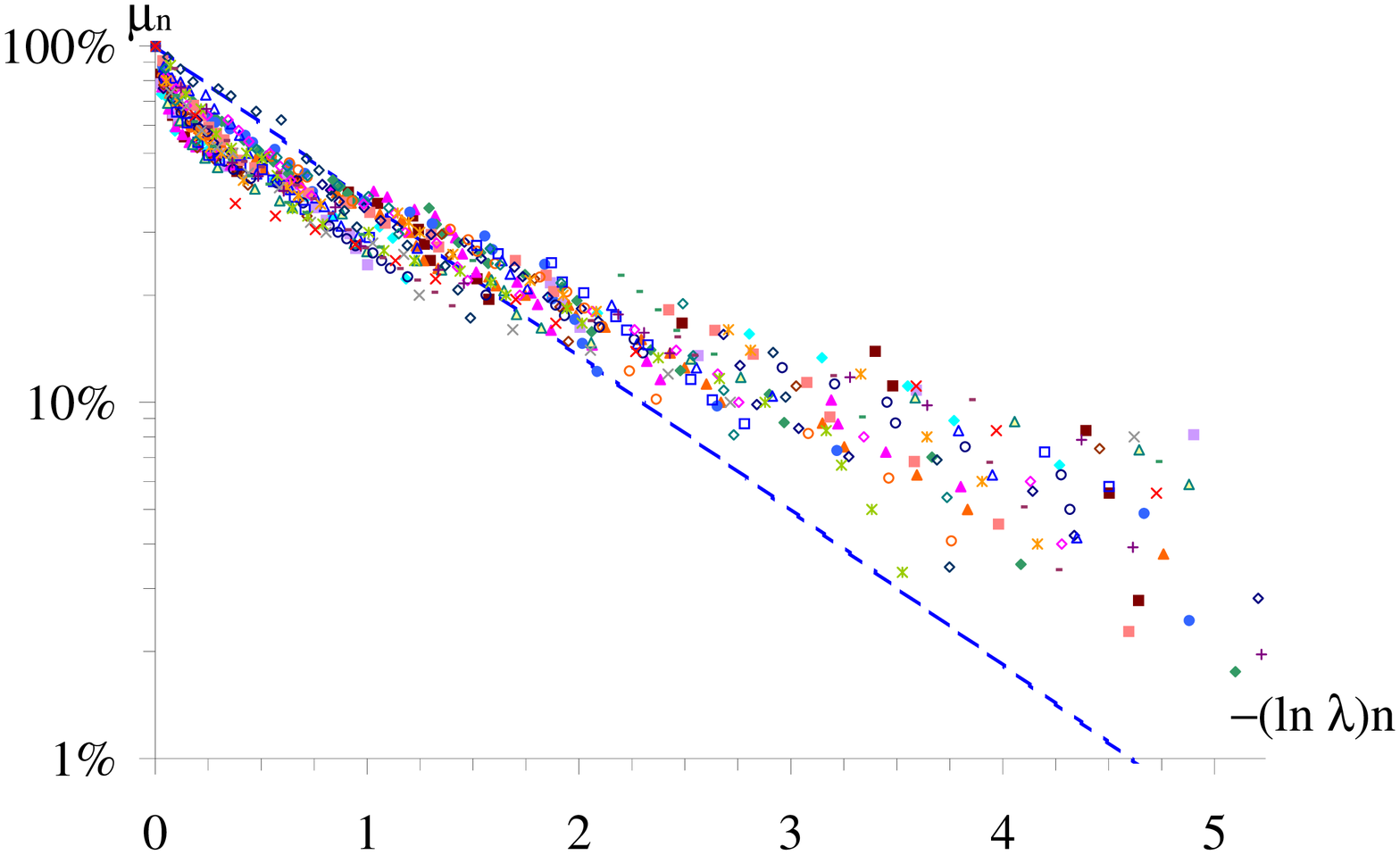}
\put(-0.9,3.5){(a)}
}
\subfigure{
    \includegraphics[height=4.1cm]{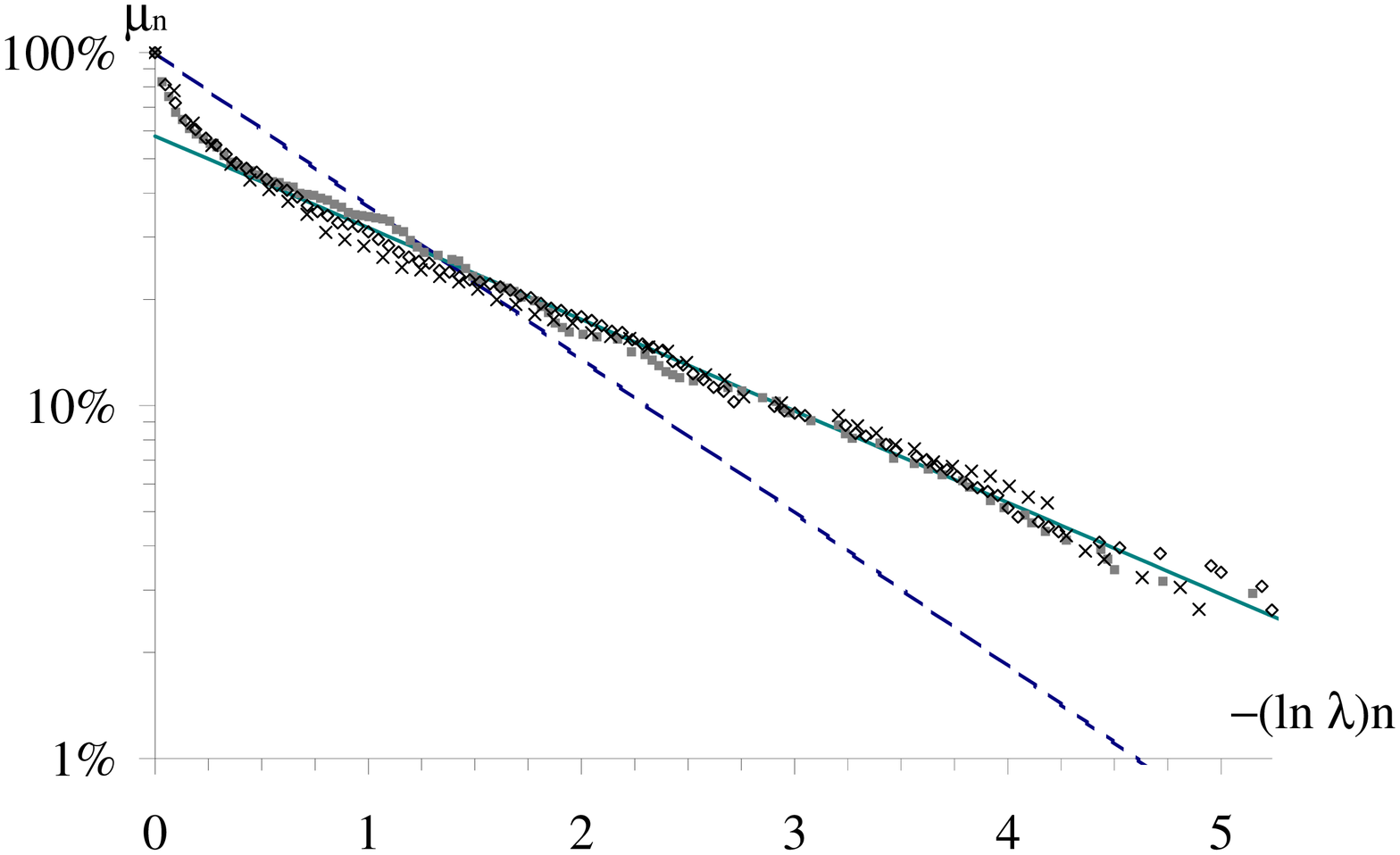}
    \label{fig-ratio(b)}
\put(-0.9,3.5){(b)}
}
\end{picture}
\caption{\label{fig-ratio}The probability to score at least $n$ points, $\mu_n$, as a function of $(-\ln \lambda)n$, with $\lambda=m/(m+1)$. (a): data for 25 French 47/2 balkline players with averages between 5 and 35. (b): three composite players (full symbols: players with an average higher than 25, empty symbols: $15 < m < 25$, crosses: $m < 15$). The dashed lines correspond to Bernoulli model. The solid line has a slope of $-0.6$.}
\end{figure}

\section{\label{position_play}Comparing players}
We would like to use the model to compare players playing at different skill levels: are the better players better in every respect or is there {\itshape one} difference which makes them better? For such a comparison one needs to cancel the effect of the difference between their averages.

\subsection{Dedimensionnalization}
In the Bernoulli model, Eq.~(\ref{mu_n_Bernoulli}), $\mu_n=\exp(n \ln \lambda)$. Letting $\nu=(-\ln \lambda)n$, $\mu_n$ is equal to $\exp(-\nu)$. If $m$ is large, $\nu \sim n/m$: the length of the runs is thus measured as a fraction of the average. This provides a way to compare players with different averages.
Figure~\ref{fig-ratio(a)} represents $\mu_n$ as a function of $(-\ln \lambda)n$. The results are very similar even for players of quite different skill levels (data are for players with averages between 5 and 35). In order to determine what differentiates the better players from the less good, we create three composite players. We use the scores of players with an average higher than 25 to create a ``player'' who played about 400 innings and who has an average of 30.4. Likewise, we use the scores of players with an average between 15 and 25 (700 innings, $m=20.5$) and of players with an average lower than 15 (500 innings, $m=10.7$).
Figure~\ref{fig-ratio(b)} represents the scores of these three ``players''.
There is no noticeable difference between them. The variations seen in Fig.~\ref{fig-ratio(a)} therefore come from  variations between players at the same skill level, not from differences between the better and less good players.

\subsection{Difference between players at different skill levels}
The slope of the asymptote to $(\ln \mu_n)$ for large values of $(-\ln \lambda)n$ is $-(\ln \rho_2)/(\ln\lambda)$. Figure~\ref{fig-ratio(b)} shows that the slope of the asymptotes is close to $-0.6$ (solid line). Thus
\begin{equation}
	\rho_2 \approx \lambda ^{0.6}
	\label{correl-rho_lambda}
\end{equation}
\noindent for all three ``players''. As most points are scored on easy shots, there logically exists a correlation between $\rho_2$ and the average (and thus $\lambda$). As difficult shots contribute little to the average, the correlation between $\rho_1$ and the average is weak and one cannot express $\rho_1$
as a function of $\lambda$ as can be done with $\rho_2$.

If the opponent always left an easy shot then the average would be close to $\rho_2/(1-\rho_2)$. Let $m_2$ the contribution to the average of runs beginning with a shot of type~2, $m_2 \approx (1-p_0)\,\lambda ^{0.6}/(1-\lambda ^{0.6})$. 
Since $m_2/m$ is a decreasing function of $\lambda$ the weaker players mostly score when their opponents leave them easy shots whereas better players also score on runs beginning with a difficult shot.

\subsection{The case of three-cushion}
Figure \ref{fig-mu-log(b)} shows $\ln \mu_n$ as a function of $n$ for Raymond Ceulemans at the 1978 three-cushion world championship, where he played 382 innings and won all his matches and the tournament. Raymond Ceulemans ---probably the best player of all times--- has more than 100 titles (world, Europe and Belgium) and he won 17 of the 18 three-cushion world championships held between 1963 and 1980. 
Even for Ceulemans at the peak of his career three-cushion is a difficult game and, unlike what was observed in Fig.~\ref{fig-mu-log(a)} for Villiers playing balkline, the results of Ceulemans are rather close to a Bernoulli process, i.e.\ to the results of a player who would not play position at all. 
A Bernoulli process seems sufficient for three-cushion billiards.

\section{\label{sect4}Calculating the scoring probabilities}

\subsection{An underdefined system}
We would like to calculate the scoring probabilities, i.e.\ the matrix $\mt{K}$\!, of a player from his scores. However, the system is underdefined: the five unknowns (the four components of $\mt{K}$ and $p_0$) lead to only three measurable parameters ($m$ and the eigenvalues $\rho_1$ and $\rho_2$). One can therefore not know $\mt{K}$ unambiguously.
All we can do is write $\mt{K}$\! as a function of $\rho_1$, $\rho_2$, $m$, and two free parameters, e.g.\ $k_{12}$ and $p_0$.
These cannot take any arbitrary values: for any $(i,\,j)$ pair, $0<k_{ij}<1$ and $k_{1j}+k_{2j}<1$. These constraints give a range of possible values for $k_{12}$ and $p_0$ and thence for $\mt{K}$\!. In spite of these constraints there is a wide range of possible values for $\mt{K}$\!: for Bernard Villiers playing 47/2 balkline one can have 
\begin{equation}
\setlength{\arraycolsep}{1.5pt}
	\mt{K} = \left(\!\begin{array}{rcrcrcr}
		40.7\% &\;& 0.0\% \\
		0.0\% &\;& 98.0\%
	\end{array} \!\right)
\text{\quad or\quad}
		\mt{K} = \left(\!\begin{array}{rcrcrcr}
		43.5\% &\;& 4.8\% \\
		32.2\% &\;& 95.2\%
	\end{array} \!\right)\negthickspace.
	\label{intervalles-Villiers-all}
\end{equation}

\subsection{A Markov process with $N$ types of shots}
Equation~(\ref{mu_Y}) shows that, if $N=2$, $(\mu_n)$ is the sum of two geometric sequences. In the general case, $(\mu_n)$ is a sum of $N$ geometric sequences. 
If we choose a value for $N$ higher than 2 the number of unknown is $N^2+N-1$ and the number of relationships between them is $N+1$. This leaves $N^2-2$ free parameters, i.e.\ two free parameters if $N=2$ but seven if $N=3$. In Fig.~\ref{fig-mu-log(a)} there does not seem to be more than two exponentials: there is no need to have $N>2$. This is fortunate as it implies that there is no need to deal with seven or more free parameters.

\subsection{\label{subsect-easy_shots}The different types of positions}
The uncertainties in expression~(\ref{intervalles-Villiers-all}) come mostly from the wide range of possible values of $p_0$. If the value of $p_0$ could be set then $\mt{K}$\! would be known more precisely. 
Setting $p_0$ defines what easy and difficult shots are: it determines how easy a shot must be to be of type~2. If $p_0$ is small then most shots are easy whereas most shots are difficult if $p_0$ is large.
We say that ``easy shots'' are shots easier than the median of the probability to score. Then by definition 50\% of the shots are easy and 50\% are difficult. 
In the remainder of this section, we consider only cases such that $p_0 = 1/2$.

For Bernard Villiers playing 47/2 balkline, if $p_0 = 1/2$ then
\begin{equation}
\setlength{\arraycolsep}{1.5pt}
	\mt{K} = \left(\!\begin{array}{rcrcrcr}
		40.9\% &\pm& 0.2\% &\quad&  1.2\% &\pm& 1.2\%\\
		11.5\% &\pm& 1.2\% &\quad& 97.8\% &\pm& 0.2\%
	\end{array} \!\right)\negthickspace,
	\label{intervalles-Villiers-12}
\end{equation}
\noindent There remain uncertainties because after setting $p_0$ there still is a free parameter. However these uncertainties are much smaller than those of Eq.~(\ref{intervalles-Villiers-all}), especially the one on $k_{21}$.

\subsection{How players can use this model}
The page \url{http://billiards.mathieu.bouville.name/biMar/} allows players to calculate their matrix $\mt{K}$ from their scores. 
The program calculates the eigenvalues by the conjugate gradient method. It then deduces the possible ranges for the components of the matrix~$\mt{K}$\! as was done with expression~(\ref{intervalles-Villiers-12}) using the scores of B.~Villiers. The user provides data in two columns, the first one is a number of points per inning and the second column is the number of occurrences of this score.
The input and the result, the probabilities to have an easy/difficult shot after an easy/difficult shot, are straightforward enough. Players need neither understanding of the model nor mathematical background to obtain this result.

\begin{figure}
\centering
\setlength{\unitlength}{1cm}
\begin{picture}(14,4.4)
\subfigure
{
    \label{avg-rho1_40-rho2_80-p10_00}
 	\ifthenelse{\boolean{color-fig}}
		{\includegraphics[height=4.4cm]{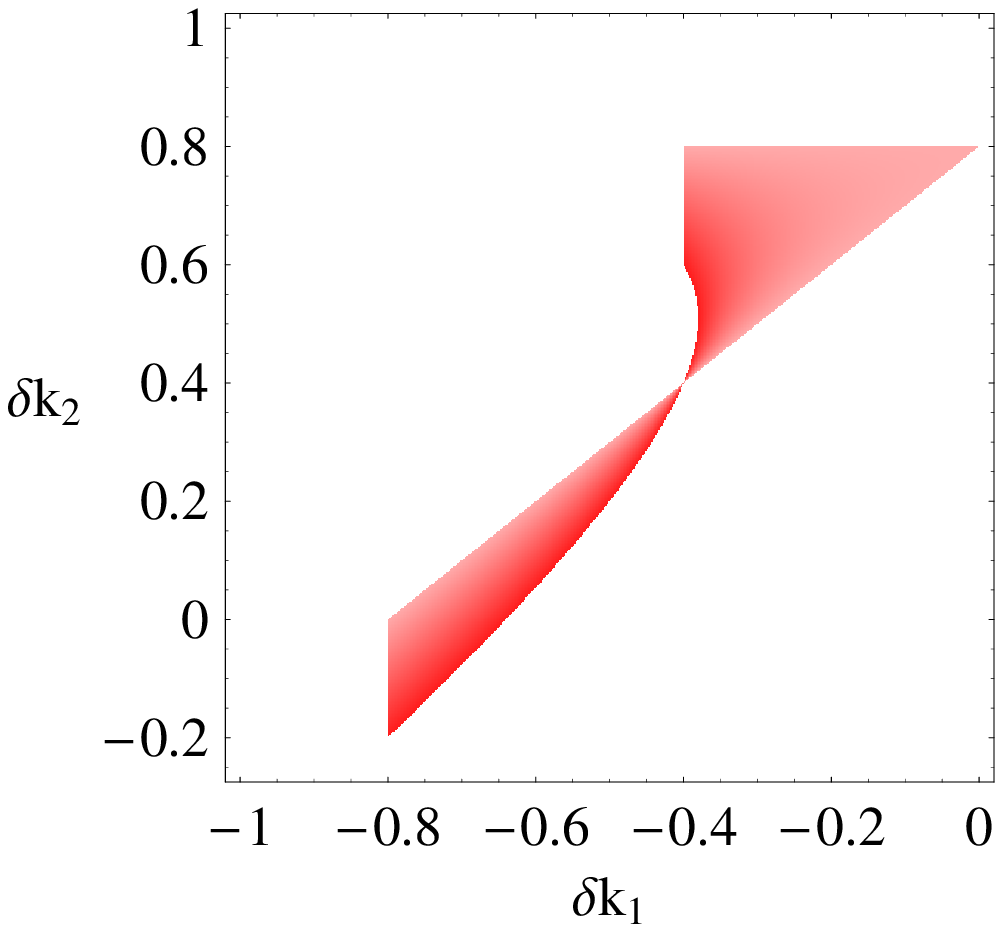}}
		{\includegraphics[height=4.4cm]{avg-rho1_40-rho2_80-p10_00}}
\put(-0.7,1){(a)}
}
\subfigure
{
    \label{avg-rho1_40-rho2_80-p10_05}
 	\ifthenelse{\boolean{color-fig}}
		{\includegraphics[height=4.4cm]{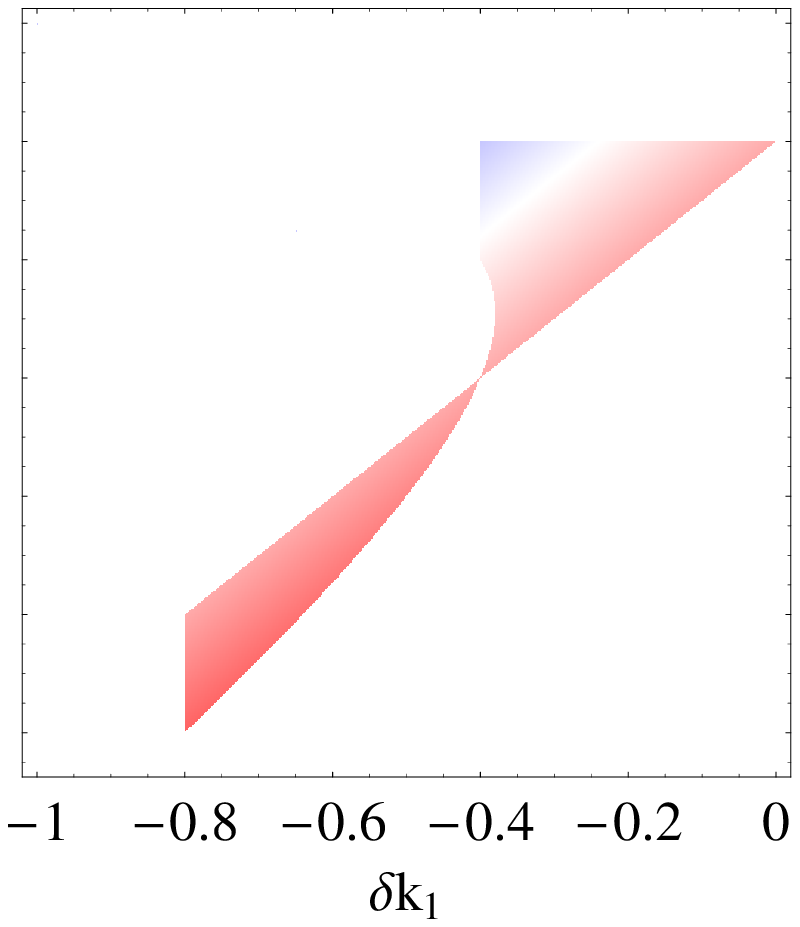}}
		{\includegraphics[height=4.4cm]{avg-rho1_40-rho2_80-p10_05}}
\put(-0.7,1){(b)}
}
\subfigure
{
    \label{avg-rho1_40-rho2_80-p10_10}
 	\ifthenelse{\boolean{color-fig}}
		{\includegraphics[height=4.4cm]{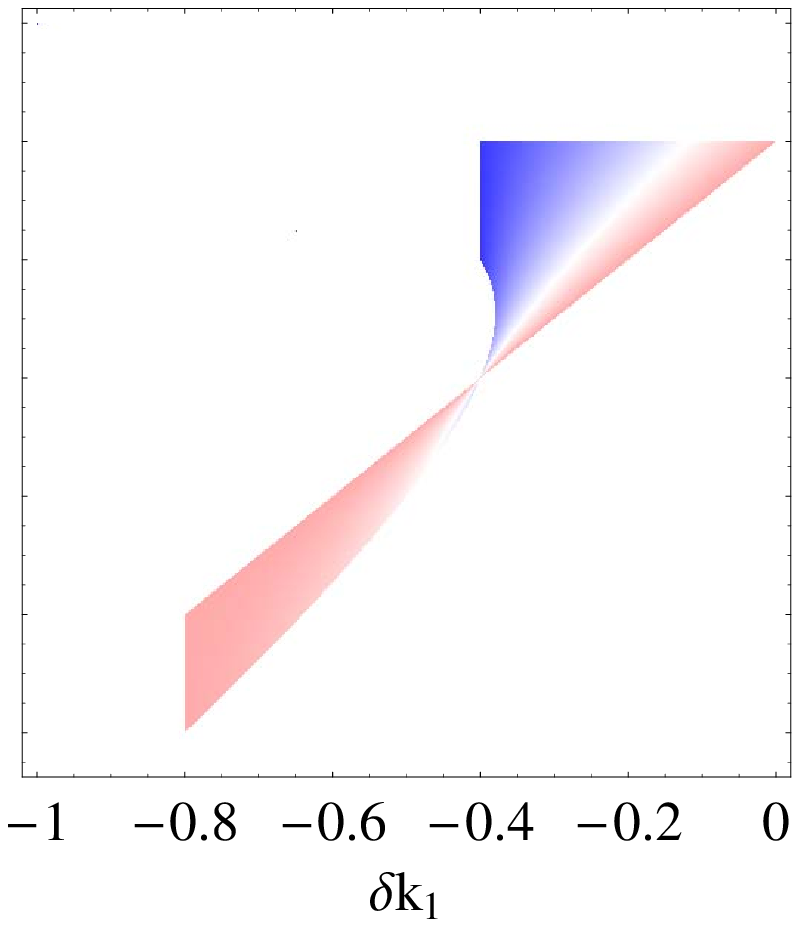}}
		{\includegraphics[height=4.4cm]{avg-rho1_40-rho2_80-p10_10}}
\put(-0.7,1){(c)}
}
 	\ifthenelse{\boolean{color-fig}}
		{\includegraphics[height=4.4cm]{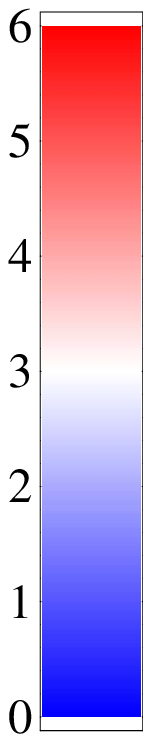}}
		{\includegraphics[height=4.4cm]{scale_6}}
\end{picture}
\ifthenelse{\boolean{color-fig}}
{\caption{\label{avg-rho1_40-rho2_80} (color) The average for $\rho_1=0.4$ and $\rho_2=0.8$. In (a) the proportion of difficult shots left by the opponent is $p_0=0$, in (b) $p_0=1/2$, and in (c) $p_0=1$.}}
{\caption{\label{avg-rho1_40-rho2_80}The average for $\rho_1=0.4$ and $\rho_2=0.8$. In (a) the proportion of difficult shots left by the opponent is $p_0=0$, in (b) $p_0=1/2$, and in (c) $p_0=1$.}}
\end{figure}

\subsection{Strategy of the opponent}
The first shot played by a player is the shot his opponent left him. It may not be random, e.g. if the opponent has a defensive style he may leave mostly difficult shots. 
Figure~\ref{avg-rho1_40-rho2_80} shows the average if $\rho_1=0.4$ and $\rho_2=0.8$, for three types of opponents. Since the average is not set there are two free parameters; we use $\delta k_1 = k_{21}-k_{11}$ and $\delta k_2 = k_{22}-k_{12}$.
In Fig.~\ref{avg-rho1_40-rho2_80-p10_00} the first shot the player has to play is always easy. In Fig.~\ref{avg-rho1_40-rho2_80-p10_10} the opponent leaves only hard shots and half easy half hard in Fig.~\ref{avg-rho1_40-rho2_80-p10_05}. If $\delta k_1=-0.4$ and $\delta k_2=0.8$ the average varies tremendously depending on the strategy of the opponent. For $\delta k_1=-0.8$ and $\delta k_2=0$ on the other hand the average is independent of the opponent. This provides extra information: if for a given player one knows $\rho_1$, $\rho_2$ and how he fares against various kinds of opponents one can reduce the span of possible values of~$\mt{K}\!$. 

\section{Conclusion}
We presented a new model of carom billiards which accounts for position play. Using a Markov process the probability to score is correlated to the difficulty of the previous shot. The probability to score at least $n$ points is then a sum of two geometrical sequences. The Bernoulli process which was used so far is a particular case of the present model with a single geometrical sequence.
Using this Markov model we established that better players score a greater proportion of their points from difficult positions left by their opponents than weaker players. The players can easily use the model through the web page \url{http://billiards.mathieu.bouville.name/biMar/}.

\section*{Acknowledgments}
I wish to thank A. Allenic for his useful comments, as well as J.-L. Frantz, J.-M. Fray, and R. Jewett who provided me with score sheets.

\appendix
\section*{\addcontentsline{toc}{section}{Glossary}Glossary}
\begin{description}
	\item[balkline:] games in which the balls cannot be kept in the same zone of the table for more than one or two shots in a row.
The three kinds of balkline are 47/2, 47/1, and 71/2.
	\item[cushion games:] games in which the cue ball must hit a certain number (one or three) of rails (cushions) before it hits the second object ball.
	\item[playing position:] trying not only to score but also to get a favorable position on the next shot, in order to have only easy shots to play and long runs (e.g., the rail nurse).
	\item[rail nurse:] in straight-rail, the player gathers the balls close to a rail in order to have a series of easy shots resulting in a very long run. 
	\item[run:] series of points without missing. Length of that series, as in ``a run of 10''.
	\item [straight-rail billiards
:] the easiest version of the game in which a point is scored when the cue ball contacts both object balls, without any additional constraint.
\end{description}

\addcontentsline{toc}{section}{References}
\bibliography{bill}
\bibliographystyle{chicago}

\end{document}


\begin{abstract}
\addcontentsline{toc}{section}{R\'esum\'e}
Le replacement est la pierre angulaire des jeux de séries au billard carambole (billard français)~: sur un coup facile le joueur peut marquer tout en s'assurant que le coup suivant sera favorable. Un coup facile est ainsi g\'en\'eralement suivi d'un autre coup facile. 
La difficult\'e d'un coup dépend donc du coup précédent. On pr\'esente un processus de Markov qui tient compte de ces corr\'elations. Ce mod\`ele permet d'expliquer les longues séries de coups faciles et les scores élevés qui en r\'esultent. 
Il permet aussi d'identifier des différences dans la nature des points sur lesquels des joueurs de niveau diff\'erent marquent le plus.
 Les joueurs peuvent utiliser ce modèle {\itshape via} la page \url{http://billiards.mathieu.bouville.name/biMar/}\!. 

\keywords{billard, sport, processus de Markov, processus de Bernoulli}
\end{abstract}

\vspace{-2ex}

\section{Introduction}
\vspace{-1ex}
Le billard carambole (ou \og billard français\fg\!) se joue avec trois billes sur une table sans poches (un court glossaire est proposé en annexe.) Si la bille du joueur touche les deux autres billes le joueur marque un point et rejoue. Quand il manque c'est au tour de son adversaire de jouer. A la fin de la partie on peut calculer une moyenne en divisant le nombre de points marqués par le nombre de reprises. Des classements peuvent ainsi être établis, même pour des joueurs qui n'ont jamais joué l'un contre l'autre; comme en athlétisme où temps et distances sont comparés pour classer des athlètes qui ne se sont jamais rencontrés.

\vspace{-1ex}

\subsection{Historique}
En 1880, à Paris, l'américain Slosson a marqué 1~103 points d'affilé face au français Vignaux qui a alors marqué 1~531 points. Dix ans plus tard Schäfer marquait 3~000 points. Ceci correspond à jouer pendant des jours sans manquer. 
De tels scores sont possibles car quand les joueurs parviennent à regrouper les billes près d'une bande et à les garder ensemble ils peuvent avoir de très longues séries. Cette technique s'appelle l'américaine. Les joueurs ne pourraient pas jouer pendant des heures sans manquer s'ils devaient faire face à de nombreux points difficiles. 
Pour empêcher les joueurs d'utiliser l'américaine ---et ainsi éviter les parties fastidieusement longues--- de nouvelles règles ont été introduites. Dans les jeux de cadres il est interdit de laisser les billes dans un certain secteur de la table pendant plus d'un ou deux coups~; de longues s\'eries (centaines de points) sont encore possible grâce au replacement. Une autre manière de compliquer le jeu est d'imposer au joueur de toucher un certain nombre de bandes avant que sa bille n'entre en contact avec la troisième bille.  

\vspace{-1ex}

\subsection{\label{sect-Bernoulli-fr}Le billard carambole comme processus de Bernoulli}
\begin{figure}
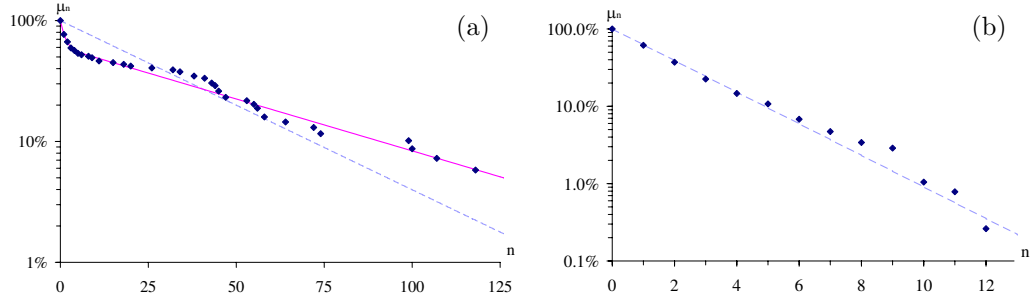

\centering
\setlength{\unitlength}{1cm}
\begin{picture}(14,4)
\subfigure
{
    \label{fig-mu-log(a)}
    \includegraphics[height=4cm]{mu-log-villiers}
\put(-0.9,3.5){(a)}
}
\subfigure
{
    \label{fig-mu-log(b)}
    \includegraphics[height=4cm]{mu-log-ceulemans}
\put(-0.9,3.5){(b)}
}
\end{picture}
\caption{\label{fig-mu-log}$\mu_n$, la probabilité de marquer au moins $n$ points, en fonction de $n$ pour Villiers au cadre 47/2 (a) et Ceulemans au trois-bandes (b). Losanges: données ; lignes pointillées : modèle bernoullien ; ligne continue: modèle markovien.}
\end{figure}

Le billard a été étudié en détail d'un point de vue mécanique par \citet{Coriolis-fr} puis plus récemment par \citet{Petit}. L'étude des statistiques du billard  a en revanche été limitée à supposer un processus de Bernoulli. 
Soit $\mu_n$ la probabilité de marquer au moins $n$ points. Pour un processus de Bernoulli avec un taux de réussite $\lambda$ le terme général de ($\mu_n$) est
\begin{equation}
	\mu_n=\lambda^n.
	\label{mu_n_Bernoulli}
\end{equation}
\noindent La longueur moyenne d'une série ---ci-apr\`es appel\'ee \og la moyenne \fg--- vaut $m=\lambda/(1-\lambda)$. La figure \ref{fig-mu-log(a)} représente des résultats de Bernard Villiers au cadre 47/2, basés sur 69 séries complètes (les séries qui ont été interrompues parce que Villiers a atteint 200 points et gagné le match ne sont pas comptabilisées). Le modèle bernoullien (pointillés) donne des résultats très différents des données (losanges)~: il surestime l'importance des courtes séries et sous-estime celle des longues séries.

\vspace{-1ex}

\subsection{L'importance du replacement}
Un processus de Bernoulli ne tient pas compte du replacement, qui est très important dans les jeux de séries. La probabilité de marquer un point dépend de sa difficulté, qui n'est pas nécessairement aléatoire.  
Sur les coups faciles, marquer n'est pas un problème et le joueur peut essayer de jouer d'une manière telle que le point sera marqué et la position finale des billes sera favorable. Ceci s'appelle le replacement. Après un coup facile, on peut donc obtenir un coup plus facile que la moyenne. 
Un joueur peut avoir une s\'erie de coups faciles et de très longues séries. De telles corrélations sont complètement absentes du modèle bernoullien décrit ci-dessus ; il ne peut donc  décrire correctement les jeux de séries.
Si quelques limitations du modèle bernoullien ont été considérées \citep{Fray-fr, Bouville-serie-fr} l'effet du replacement n'a jamais été étudié.

\vspace{-1ex}

\section{\label{sect-Markov-fr}Le billard carambole comme processus de Markov}
\subsection{D\'efinir un processus de Markov}
Si on considère qu'il existe $N_0$ positions possibles des billes sur la table, le résultat d'un coup peut être un succès avec les billes finissant dans une des $N_0$ configurations possibles ou un échec. La probabilité de marquer et la position des billes pour le prochain coup dépendent seulement de la position des billes au début du coup, pas de leur position antérieure. La succession des coups est donc un processus de Markov à $N_0+1$ états
. Malheureusement les probabilités de transition sont nombreuses, $O(N_0{}^2)$, et ne peuvent être aisément obtenues. 
Pour simplifier le mod\`ele et le rendre utilisable, on peut regrouper des positions similaires pour n'avoir que $N$ \emph{types} de positions. Si $N$ est petit, les paramètres sont suffisamment peu nombreux pour pouvoir être calculés et interprétés concrètement. 

\vspace{-1ex}

\subsection{La probabilit\'e de marquer au moins $n$ points}
Soit $\vc{P}_{\!n}$ le vecteur \`a $N$ dimensions dont la composante $i$ est la probabilité que les billes soient dans une position de type $i$ après le coup $n$ (on utilise \og coup de type $i$ \fg\ et \og position [des billes sur la table] de type $i$ \fg\ comme des synonymes). Pour tout $n\ge0$ la suite $\left(\vc{P}_{\!n} \right)$ obéit à $\vc{P}_{\!n+1} = \mt{K} \vc{P}_{\!n}$. L'élément $k_{ij}$ de $\mt{K}$\!, matrice $N \times N$\!, est la probabilité de marquer sur un coup de type $j$ et que le prochain coup soit de type $i$. $\mt{K}$\! est {\itshape apparentée} à une matrice de transitions~: l'état correspondant à un manqué n'est pas inclus car les composantes afférentes peuvent être obtenues simplement à partir de $\mt{K}$\!. La probabilité de manquer un coup de type $i$ est $1-\kappa_i$, o\`u $\kappa_i=\sum_j k_{ji}$ est la probabilité de marquer sur un coup de type~$i$. 

Dans ce qui suit on utilise $N=2$ mais le mod\`ele peut ais\'ement \^etre \'etendu \`a plus de types de coups. Soient $\rho_1<\rho_2$ les valeurs propres de $\mt{K}$\! et $\vc{V}_{\!1}$ et $\vc{V}_{\!2}$ les vecteurs propres associés. Pour tout $n \ge 0$, $\vc{P}_{\!n}$ s'écrit
\begin{equation}
	\vc{P}_{\!n} = \left(\begin{array}{cc}	\vc{V}_{\!1} & \vc{V}_{\!2} \end{array} \right)
		\left(\begin{array}{cc}
		\rho_1{}^n & 0\\
		 0 & \rho_2{}^n
	\end{array} \right)
 		\left(\begin{array}{cc}	\vc{V}_{\!1} & \vc{V}_{\!2} \end{array} \right)^{-1}
		\left(\begin{array}{cc}	p_0 \\ 1-p_0 \end{array} \right)\negthickspace.
	\label{eq_P}
\end{equation}
\noindent $p_0$ est la probabilité que le premier point joué par le joueur soit difficile. 
La probabilité de marquer au moins $n$ points, $\mu_n$, est la somme des composantes de $\vc{P}_{\!n}$. Pour tout $n \ge 0$, $\mu_n$ est de la forme 
\begin{equation}
	\mu_n = (1-Y) \rho_1^n + Y \rho_2^n,
	\label{mu_Y}
\end{equation}
\noindent $Y$ est un scalaire connu en fonction de $\vc{V}_{\!1}$, $\vc{V}_{\!2}$ et $p_0$. La figure~\ref{fig-mu-log(a)} présente $\mu_n$ sur une échelle semi-logarithmique en fonction de $n$ pour Bernard Villiers au cadre 47/2~: les données suivent une exponentielle pour les petits $n$ et une autre aux valeurs élevées de $n$. Le processus de Markov, équation~(\ref{mu_Y}), est une meilleure description des jeux de cadres que le processus de Bernoulli, équation (\ref{mu_n_Bernoulli}).

La moyenne, $m$, est la somme sur $n$ de $\mu_n$,
\begin{equation}
	m = (1-Y) \frac{\rho_1}{1-\rho_1} + Y \frac{\rho_2}{1-\rho_2}.
	\label{def-m}
\end{equation}
\noindent Comme $Y$ est connu en fonction de $\rho_1$, $\rho_2$ et $m$, l'\'equation~(\ref{mu_Y}) peut \^etre r\'e\'ecrite en exprimant $\mu_n$ en fonction de $\rho_1$, $\rho_2$, $m$ et $n$.
Comme la moyenne est connue grâce aux feuilles de matches $(\mu_n)$ ne dépend que de deux inconnues, $\rho_1$ et $\rho_2$. On trouve leur valeur en minimisant la somme des carr\'es des diff\'erences entre la valeur de $\mu_n$ des donn\'ees et celle de l'\'equation~(\ref{mu_Y}).

\vspace{-1ex}

\subsection{Types de coups}
On a consid\'er\'e qu'il existe deux types de coups, d\'efini la matrice~$\mt{K}$\! et utilis\'e ses vecteurs propres $\rho_1$ et $\rho_2$ pour calculer $\mu_n$ dans l'\'equation~(\ref{mu_Y}). Mais on n'a jamais vraiment sp\'ecifi\'e ce que sont ces deux types de coups. Ils ne peuvent pas \^etre al\'eatoires, sinon le processus de Markov se r\'eduirait \`a un processus de Bernoulli. Les positions de type 1 (resp.\ 2) doivent \^etre telles qu'elles sont g\'en\'eralement suivies d'une autre position de type 1 (resp.\ 2).
A la libre, les coups de type 2 seraient des positions d'am\'ericaine: il y a une grande probabilit\'e qu'une position d'am\'ericaines soit suivie d'une autre position d'am\'ericaine. Les coups de type 1 seraient les positions moins favorable des billes sur la table. 
Plus g\'en\'eralement les positions de type 1 sont \og difficiles \fg\ et celles de type 2 \og faciles \fg\!.

\begin{figure}
\centering
\setlength{\unitlength}{1cm}
\begin{picture}(14,4.6)
\subfigure{
    \label{fig-ratio(a)}
    \includegraphics[height=4.1cm]{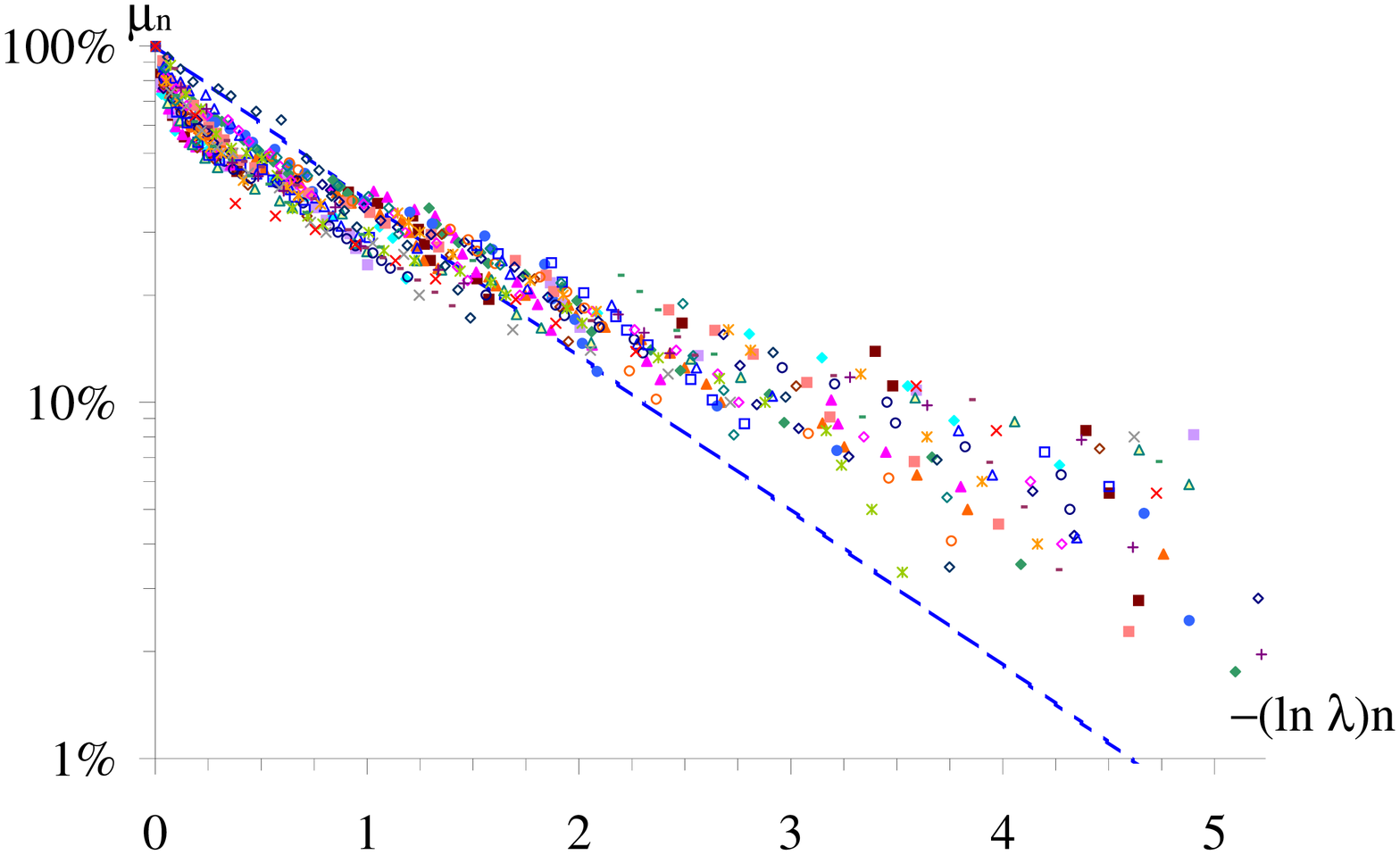}
\put(-0.9,3.5){(a)}
}
\subfigure{
    \includegraphics[height=4.1cm]{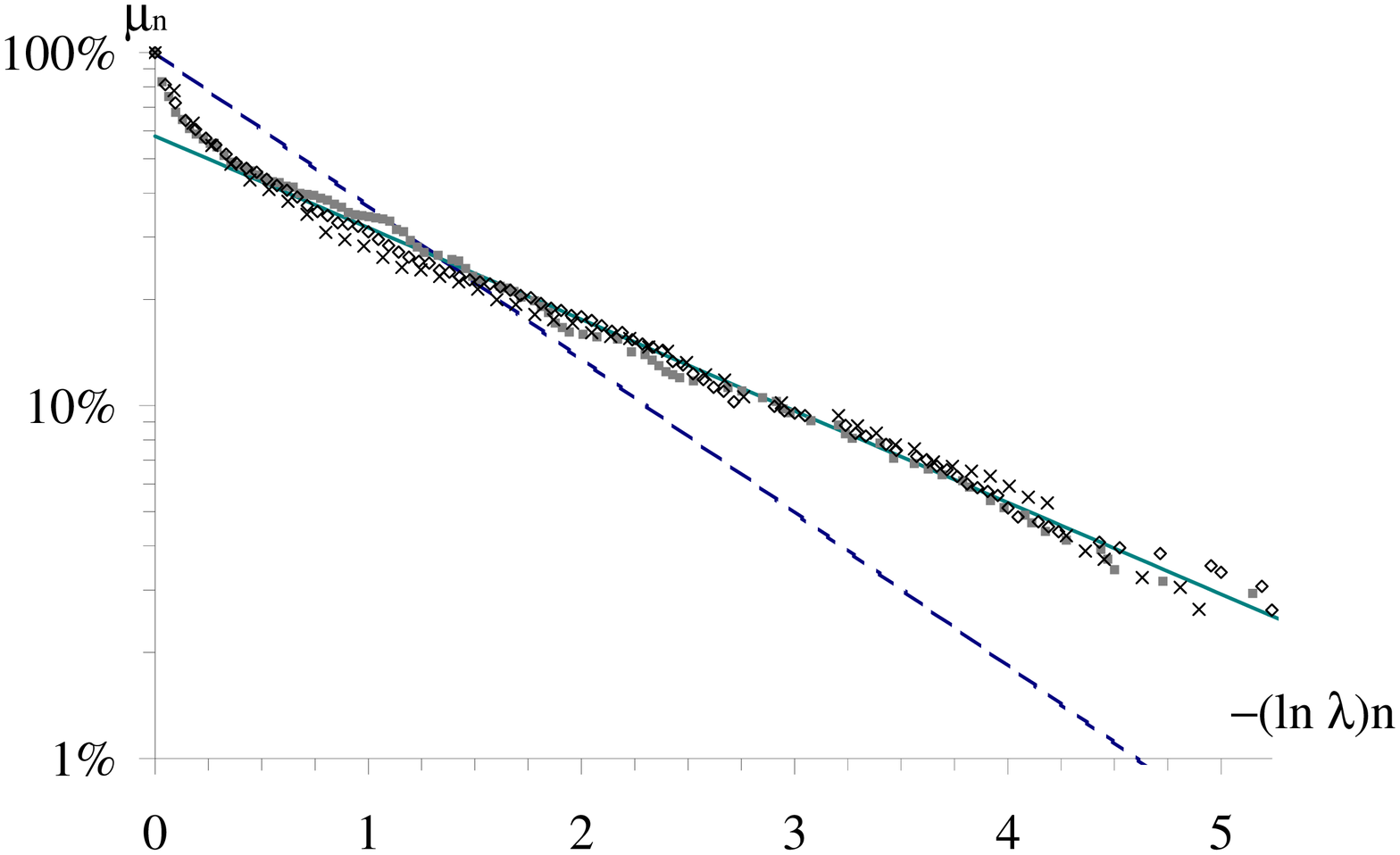}
    \label{fig-ratio(b)}
\put(-0.9,3.5){(b)}
}
\end{picture}
\caption{\label{fig-ratio}La probabilit\'e de marquer au moins $n$ points, $\mu_n$, en fonction de $(-\ln \lambda)n$. (a) : donn\'ees de 25 joueurs français de cadre 47/2. (b) : trois joueurs composites (symboles pleins: joueurs de moyenne sup\'erieure \`a 25, symboles vides : $15 < m < 25$, croix: $m < 15$). Les lignes pointill\'ees correspondent au mod\`ele bernoullien. La ligne continue a une pente de $-0,6$.}
\end{figure}

\vspace{-1ex}

\section{\label{replacement}Comparaisons de joueurs}
\vspace{-1ex}
On aimerait utiliser le mod\`ele afin de comparer des joueurs de niveaux diff\'erents : les meilleurs sont-ils meilleurs en out ou bien y a-t-il {\itshape une} difference qui les rend meilleurs ? Pour une telle comparaison il faut annuler l'effet de la diff\'erence de leurs moyennes.

\vspace{-1ex}

\subsection{D\'edimensionnalisation}
Dans le mod\`ele bernoullien, \'equation~(\ref{mu_n_Bernoulli}), $\mu_n=\exp(n \ln \lambda)$. En posant $\nu=(-\ln \lambda)n$,  $\mu_n$ est \'egal \`a $\exp(-\nu)$. Si $m$ est grand, 
$\nu \sim n/m$~; la longueur des s\'eries est ainsi mesur\'ee en pourcentage de la moyenne.
La figure~\ref{fig-ratio(a)} repr\'esente $\mu_n$ en fonction de $(-\ln \lambda)n$. Les r\'esultats sont tr\`es similaires m\^eme pour des joueurs de niveaux tr\`es diff\'erents (les donn\'ees sont pour des joueurs de moyenne comprise entre 5 et 35). 
Afin de d\'eterminer ce qui diff\'erencie meilleurs et moins bons, on cr\'ee trois joueurs composites. On utilise les scores de joueurs de moyenne sup\'erieure \`a 25 afin de cr\'eer un \og joueur~\fg\ ayant jou\'e environ 400 reprises pour une moyenne de 30,4. De m\^eme on utilise les scores de joueurs ayant entre 15 et 25 de moyenne (700 reprises, $m=20,5$) et de joueurs de moyenne inf\'erieure \`a 15 (500 reprises, $m=10,7$).
La figure~\ref{fig-ratio(b)} repr\'esente les scores de ces trois \og joueurs~\fg\!. 
Il n'y a pas de diff\'erence notable entre eux. Les variations de la figure~\ref{fig-ratio(a)} proviennent donc de variations entre joueurs de m\^eme niveau et non de diff\'erences entre meilleurs et moins bons.

\vspace{-1ex}

\subsection{Diff\'erence entre joueurs de niveau diff\'erent}
La pente de l'asymptote de $(\ln \mu_n)$ pour les grandes valeurs de $(-\ln \lambda)n$ est $-(\ln \rho_2)/(\ln\lambda)$. La figure~\ref{fig-ratio(b)} montre que la pente des asymptotes est $-0,6$ (ligne continue). On a donc
\begin{equation}
	\rho_2 \approx \lambda ^{0,6}
	\label{correl-rho_lambda}
\end{equation}
\noindent pour ces trois \og joueurs~\fg\!. Comme la plupart des points sont marqu\'es sur des coups faciles, il existe logiquement une corr\'elation entre $\rho_2$ et la moyenne, donc $\lambda$. Les coups difficiles contribuant peu \`a la moyenne, la corr\'elation entre $\rho_1$ et la moyenne est faible et on ne peut exprimer $\rho_1$ en fonction de $\lambda$ comme on peut le faire avec $\rho_2$.

Si l'adversaire laissait toujours un coup facile la moyenne serait d'environ $\rho_2/(1-\rho_2)$. Soit $m_2$ la contribution \`a la moyenne des s\'eries commençant par un coup de type 2, $m_2 \approx (1-p_0)\,\lambda ^{0,6}/(1-\lambda ^{0,6})$.
Comme $m_2/m$ est une fonction d\'ecroissante de $\lambda$  les joueurs les moins bons marquent principalement quand leur adversaire leur laisse un coup facile alors que les meilleurs  marquent aussi sur des s\'eries faites \`a partir d'une position difficile.

\vspace{-1ex}

\subsection{Cas du trois-bandes}
La figure \ref{fig-mu-log(b)} représente $\ln \mu_n$ en fonction de $n$ pour Raymond Ceulemans au championnat du monde de trois-bandes 1978 o\`u il a joué 382 reprises, remportant tous ses matches et le titre. Raymond Ceulemans ---probablement le meilleur joueur de tous les temps--- a plus de 100 titres (mondiaux, européens et belges) et a gagné 17 championnats du monde de trois-bandes sur 18 entre 1963 et 1980. 
Même pour Ceulemans au sommet de sa carrière le trois-bandes est difficile et, à la différence de ce qui a été observé sur la figure~\ref{fig-mu-log(a)} pour Villiers jouant au cadre, les r\'esultats de Ceulemans sont proches d'un processus de Bernoulli, c'est-\`a-dire des r\'esultats d'un joueur qui ne jouerait pas du tout le replacement. 
Le modèle bernoullien semble donc suffire dans le cas du trois-bandes.

\vspace{-1ex}

\section{\label{sect4}Calcul des probabilit\'es de marquer}
\subsection{Un système sous-défini}
On aimerait calculer les probabilit\'es de marquer (la matrice $\mt{K}$) \`a partir de ses scores. Cependant le système est sous-défini : les cinq inconnues (les quatre éléments de $\mt{K}$\! et $p_0$) ne débouchent que sur trois param\`etres mesurables ($m$ et les valeurs propres $\rho_1$ et $\rho_2$). On ne peut donc pas dé\-ter\-miner $\mt{K}$ uniquement.
On ne peut qu'\'ecrire $\mt{K}$ en fonction de $\rho_1$, $\rho_2$, $m$ et de deux param\`etres libres, par exemple $k_{12}$ et $p_0$.
Ils ne peuvent pas prendre n'importe quelle valeur : pour tout  $(i,\,j)$ on doit avoir $0<k_{ij}<1$ et $k_{1j}+k_{2j}<1$. 
Ces contraintes donnent un intervalle de valeurs possibles pour $k_{12}$ et $p_0$ puis pour $\mt{K}$\!.
Malgr\'e ces contraintes il existe de nombreuses valeurs possibles de $\mt{K}$\!: pour B.~Villiers au cadre 47/2 on peut avoir
\begin{equation}
\setlength{\arraycolsep}{1.5pt}
	\mt{K} = \left(\!\begin{array}{rcrcrcr}
		40,7~\% &\;& 0,0~\% \\
		0,0~\% &\;& 98,0~\%
	\end{array} \!\right)
\text{\quad ou\quad}
		\mt{K} = \left(\!\begin{array}{rcrcrcr}
		43,5~\% &\;& 4,8~\% \\
		32,2~\% &\;& 95,2~\%
	\end{array} \!\right)\negthickspace.
	\label{intervalles-Villiers-all}
\end{equation}

\vspace{-1ex}

\subsection{Un processus de Markov avec $N$ types de coups}
L'\'equation~(\ref{mu_Y}) montre que si $N=2$ $(\mu_n)$ est la somme de deux suites g\'eom\'etriques. Dans le cas g\'en\'eral, $(\mu_n)$ est la somme de $N$ suites g\'eom\'etriques. 
Si on choisit une valeur de $N$ sup\'erieure \`a 2 le nombre d'inconnues est de $N^2+N-1$ et le nombre de relations entre elles de $N+1$. Ceci laisse $N^2-2$ param\`etres libres~: deux param\`etres libres si $N=2$ mais sept si $N=3$. Sur la figure~\ref{fig-mu-log(a)} il ne semble pas y avoir plus de deux exponentielles~: on n'a pas besoin d'avoir $N>2$. C'est heureux car \c{c}a signifie qu'il n'est pas n\'ecessaire d'avoir sept param\`etres libres (ou plus).

\vspace{-1ex}

\subsection{\label{subsect-easy_shots}Les diff\'erents types de positions}
Les incertitudes dans l'expression~(\ref{intervalles-Villiers-all}) proviennent surtout du large intervalle de valeurs possibles pour $p_0$. Si on fixe $p_0$ alors $\mt{K}$\! est connue plus pr\'ecis\'ement. 
La valeur de $p_0$ d\'efinit ce que sont les coups faciles et difficiles : elle determine quelle niveau de facilit\'e est requis pour qu'une position soit de type 2. Si $p_0$ est petit la plupart des coups sont faciles alors que la plupart des coups sont difficiles si $p_0$ est grand.
On dit que les \og coups faciles\fg\ sont des coups plus faciles que la médiane de la probabilité de marquer. Alors par définition 50~\% des coups sont faciles et 50~\% difficiles. 
Dans ce qui suit on se limite donc au cas $p_0 = 1/2$.

Pour Bernard Villiers jouant au cadre 47/2, si $p_0=1/2$
\begin{equation}
\setlength{\arraycolsep}{1.5pt}
	\mt{K} = \left(\!\begin{array}{rcrcrcr}
		40,9~\% &\pm& 0,2~\% &\;& 1,2~\% &\pm& 1,2~\%\\
		11,5~\% &\pm& 1,2~\% &\;& 97,8~\% &\pm& 0,2~\%
	\end{array} \!\right)\negthickspace.
	\label{intervalles-Villiers-12-min}
\end{equation}
\noindent Des incertitudes demeurent vu qu'apr\`es avoir fix\'e $p_0$ il y a encore un param\`etre libre. Ces incertitudes sont tout de m\^eme bien plus faibles que celles de l'\'equation~(\ref{intervalles-Villiers-all}), surtout $k_{21}$.

\vspace{-1ex}

\subsection{Comment les joueurs peuvent utiliser ce modèle}
La page \url{http://billiards.mathieu.bouville.name/biMar/}\! permet aux joueurs de calculer leur matrice $\mt{K}$ à partir de leurs scores. 
Le programme calcule les valeurs propres par la méthode du gradient conjugué. Il en déduit l'intervalle possible des \'el\'ements de la matrice $\mt{K}$\! comme on a fait pour B.~Villiers dans l'expression~(\ref{intervalles-Villiers-12-min}). L'utilisateur fournit des données en deux colonnes, la première est un nombre de points par reprise et la seconde est le nombre de fois où ce score a été obtenu.
Ces données ainsi que le résultat, la probabilité d'obtenir un coup facile/difficile après un coup facile/difficile, sont assez simples. Le joueur n'a pas besoin de conna\^itre le modèle ni les mathématiques pour obtenir ce résultat.

\begin{figure}
\centering
\setlength{\unitlength}{1cm}
\begin{picture}(14,4.4)
\subfigure
{
    \label{avg-rho1_40-rho2_80-p10_00}
 	\ifthenelse{\boolean{color-fig}}
		{\includegraphics[height=4.4cm]{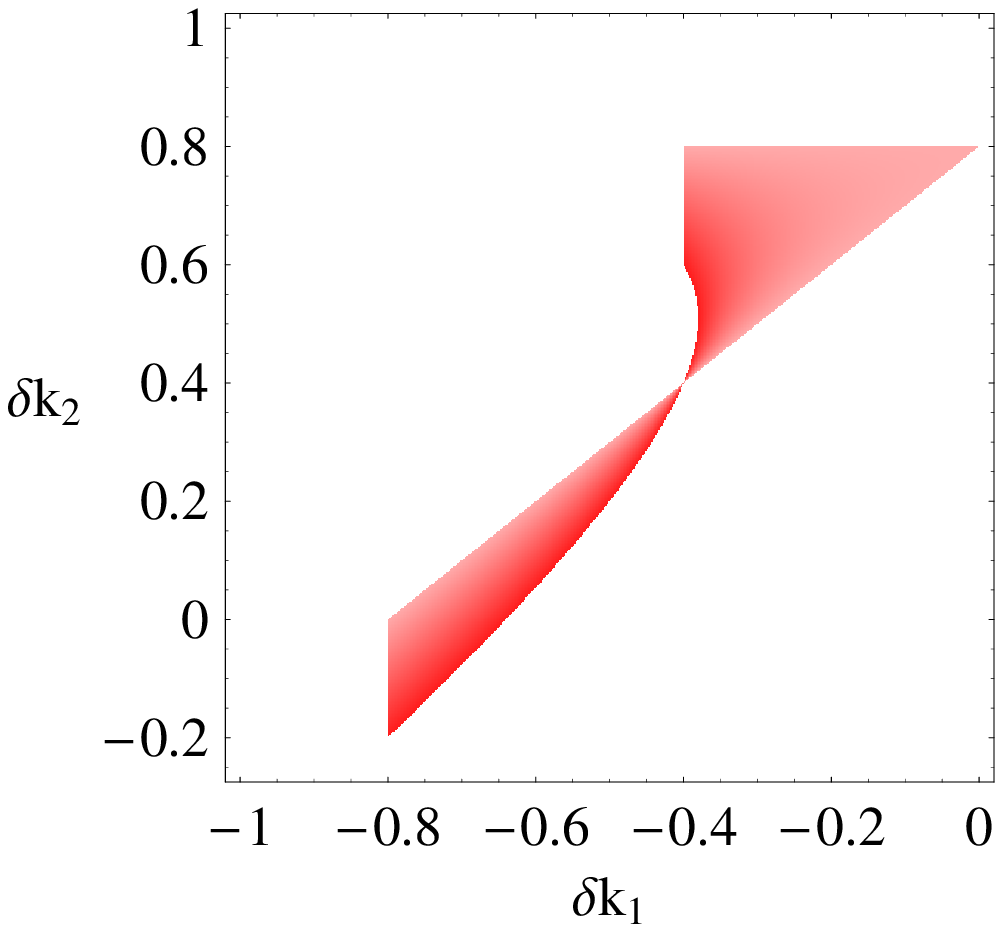}}
		{\includegraphics[height=4.4cm]{avg-rho1_40-rho2_80-p10_00}}
\put(-0.7,1){(a)}
}
\subfigure
{
    \label{avg-rho1_40-rho2_80-p10_05}
 	\ifthenelse{\boolean{color-fig}}
		{\includegraphics[height=4.4cm]{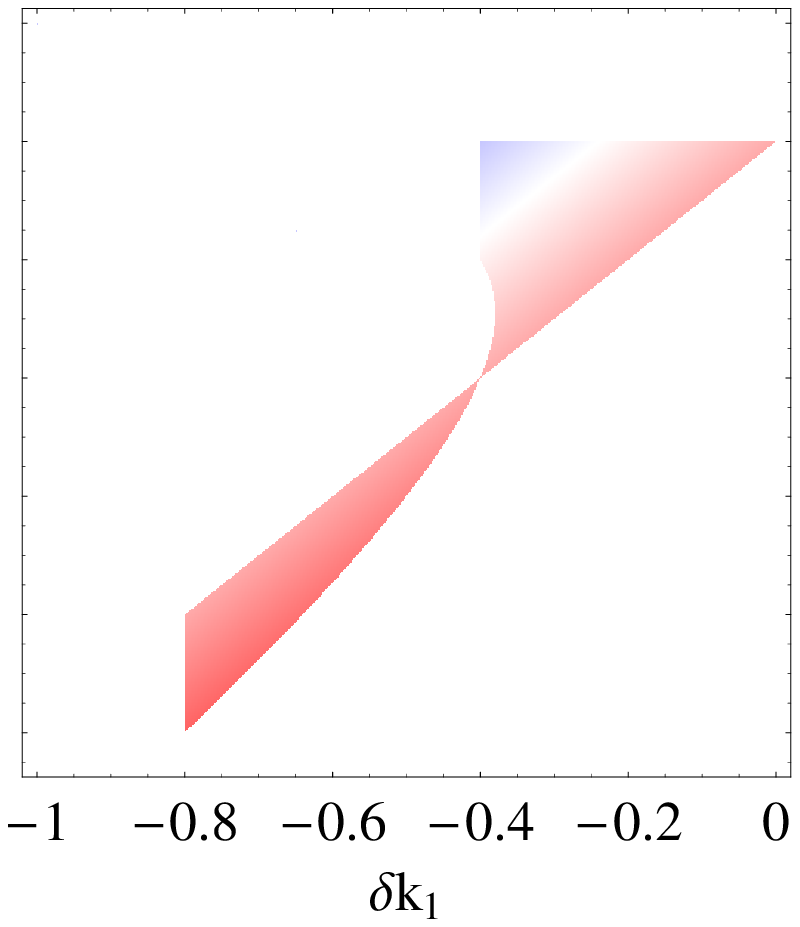}}
		{\includegraphics[height=4.4cm]{avg-rho1_40-rho2_80-p10_05}}
\put(-0.7,1){(b)}
}
\subfigure
{
    \label{avg-rho1_40-rho2_80-p10_10}
 	\ifthenelse{\boolean{color-fig}}
		{\includegraphics[height=4.4cm]{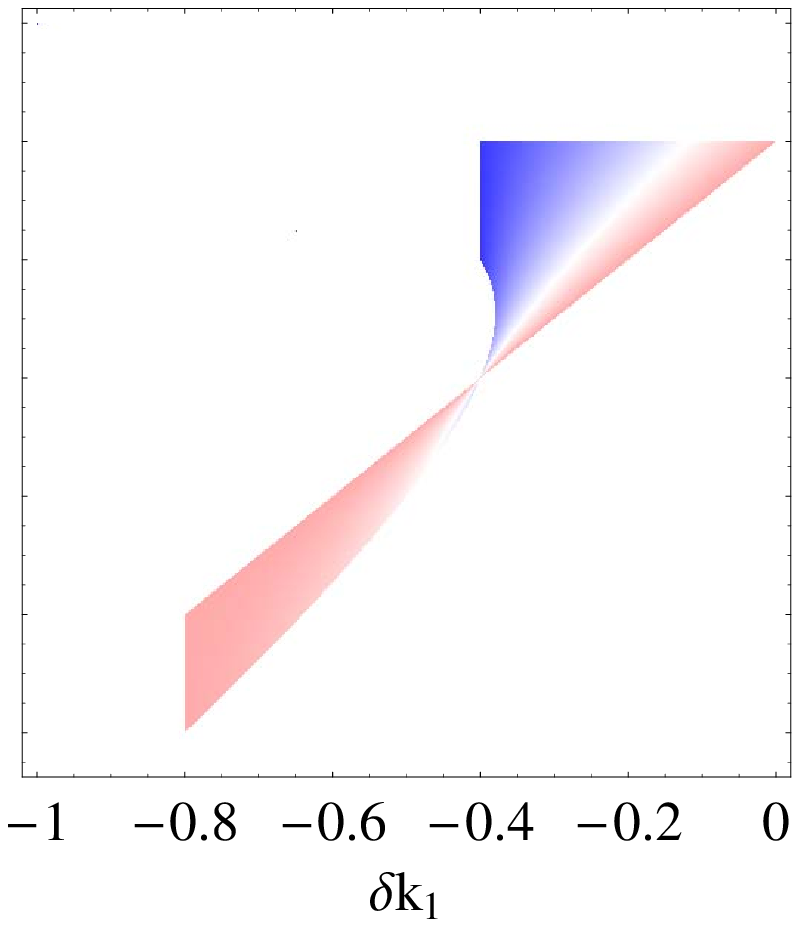}}
		{\includegraphics[height=4.4cm]{avg-rho1_40-rho2_80-p10_10}}
\put(-0.7,1){(c)}
}
 	\ifthenelse{\boolean{color-fig}}
		{\includegraphics[height=4.4cm]{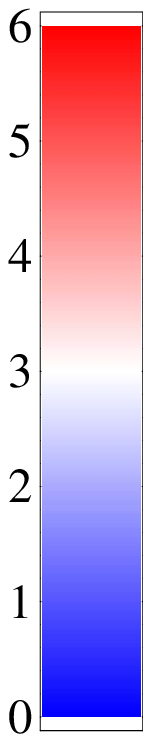}}
		{\includegraphics[height=4.4cm]{scale_6}}
\end{picture}
\ifthenelse{\boolean{color-fig}}
{\caption{\label{avg-rho1_40-rho2_80} (couleur) La moyenne pour $\rho_1=0.4$ et $\rho_2=0.8$. En (a) la proportion de coups difficiles laissés par l'adversaire est de $p_0=0$, en (b) $p_0=1/2$ et en (c) $p_0=1$.}
}
{\caption{\label{avg-rho1_40-rho2_80}La moyenne pour $\rho_1=0.4$ et $\rho_2=0.8$. En (a) la proportion de coups difficiles laissés par l'adversaire est de $p_0=0$, en (b) $p_0=1/2$ et en (c) $p_0=1$.}
}
\end{figure}

\vspace{-1ex}

\subsection{Stratégie de l'adversaire}
Le premier coup joué par un joueur est le coup laissé par son adversaire. Il n'est pas forcément aléatoire : si l'adversaire joue défensivement, il laissera souvent les billes dans une position difficile. La figure \ref{avg-rho1_40-rho2_80} représente la moyenne si $\rho_1=0,4$ et $\rho_2=0,8$, pour trois types d'adversaires. Comme la moyenne n'est pas fix\'ee il y a deux param\`etres libres~; on utilise $\delta k_1 = k_{21}-k_{11}$ et $\delta k_2 = k_{22}-k_{12}$.
 Pour la figure \ref{avg-rho1_40-rho2_80-p10_00} le premier coup que le joueur doit jouer est toujours facile.  Dans le cas de la figure \ref{avg-rho1_40-rho2_80-p10_10} l'adversaire laisse seulement des coups difficiles et moitié moitié pour la figure \ref{avg-rho1_40-rho2_80-p10_05}. Si $\delta k_1=-0,4$ et $\delta k_2=0,8$ la moyenne d\'epend fortement de la stratégie de l'adversaire.  Quand $\delta k_1=-0,8$ et $\delta k_2=0$ en revanche la moyenne est indépendante de ce que fait l'adversaire. Ceci fournit des informations supplémentaires : si pour un joueur donné on connaît $\rho_1$ et $\rho_2$ et comment il se comporte contre divers joueurs il est possible de réduire le nombre de valeurs possibles de $\mt{K}$\!. 

\vspace{-1ex}

\section{Conclusion}
\vspace{-1ex}
On a présenté un nouveau modèle du billard carambole qui tient compte du replacement. En modélisant les jeux de séries comme un processus de Markov la probabilité de marquer un point est corrélée à la difficulté du coup précédent. La probabilité de marquer au moins $n$ points est alors une somme de deux suites géométriques. Le processus de Bernoulli est un cas particulier du présent modèle avec une seule suite géométrique. En utilisant ce modèle markovien on a pu établir des différences entre joueurs de diff\'erents niveaux~: les meilleurs joueurs marquent une plus grande proportion de leurs points sur des positions difficiles laiss\'ees par leur adversaire que les moins bons. La page \url{http://billiards.mathieu.bouville.name/biMar/}\! permet aux joueurs d'utiliser ce modèle sans difficulté.

\vspace{-1ex}

\section*{Remerciements}
\vspace{-1ex}
Je tiens à remercier A. Allenic pour ses conseils avisés ainsi que J.-L. Frantz, J.-M. Fray et R. Jewett pour m'avoir fourni des feuilles de matches.

\vspace{-1ex}

\appendix
\section*{\addcontentsline{toc}{section}{Glossaire}Glossaire}
\vspace{-1ex}
\begin{description}
	\item[américaine :] à la libre, le joueur rassemble les billes à proximité d'une bande afin de n'avoir que des coups faciles à jouer et d'obtenir de très longues séries.
	\item[bandes (une bande, trois bandes) :] modes de jeu où la bille du joueur doit entrer en contact avec un certain nombre de bandes (une ou trois) avant de conclure le point.
	\item[cadres :] trois modes de jeu (47/2, 47/1 et 71/2) où les billes ne peuvent pas être laissées dans la même partie du billard pendant plus d'un ou deux coups.
	\item[partie libre (ou libre) :] mode de jeu le plus simple dans lequel un point est marqué si la bille du joueur touche les deux autres billes, sans autre contrainte.
	\item[replacement :] fait d'essayer non seulement de marquer le point mais aussi d'obtenir une position favorable pour le prochain coup afin de n'avoir que des coups faciles à jouer et des scores élevés. L'américaine en est un exemple.
	\item[série :] coups successifs sans manquer. Nombre de points ainsi marqués, {\itshape une série de 10}.
\end{description}

\vspace{-1ex}

\addcontentsline{toc}{section}{Bibliographie}
\bibliography{bill}
\bibliographystyle{chicago}